\documentclass[fleqn,10pt]{olplainarticle}
\usepackage{graphicx}
\usepackage{amsmath,amsthm}
\usepackage{soul} 

\newtheorem{lem}{Lemma}
\newtheorem{prop}{Proposition}
\newtheorem{cor}{Corollary}
\newtheorem{thm}{Theorem}

\theoremstyle{definition}
\newtheorem{ass}{Assumption}
\newtheorem{defn}{Definition}
\newtheorem{prob}{Problem}

\theoremstyle{remark}
\newtheorem{rem}{Remark}

\newcommand{\norm}[1]{\left\|#1\right\|}

\DeclareMathOperator*{\esssup}{ess\,sup}

\newcommand{\eu}{\mathrm{e}}
\newcommand{\RE}{\mathbb{R}}
\newcommand{\REp}{\overline{\RE}_+}
\newcommand{\rc}{\mathrm{c}}
\newcommand{\rd}{\mathrm{d}}

\newcommand{\CO}{\mathcal{C}}
\newcommand{\COh}{\CO([-h,0],\RE^n)}
\newcommand{\LE}{\mathcal{L}}
\newcommand{\CK}{\mathcal{K}}
\newcommand{\CKL}{\mathcal{K}\!\!\mathcal{L}}

\title{Dominant-Pole Placement for \\ Predictor Synthesis}

\author[1]{Bryan Rojas-Ricca}
\author[1]{Fernando Casta{\~n}os}
\author[1]{Sabine Mondi{\'e}}
\affil[1]{Departamento de Control Autom\'atico, Cinvestav-IPN, Mexico}


\begin{abstract}
This article analyzes the high-gain prediction approach for nonlinear input-delay systems.
The problem is discussed in the light of weighted homogeneity and input-to-state stability.
The canonical form for uniformly observable nonlinear systems allows tuning the spectrum of
the linear part by multiplicity-induced dominance and ensures closed-loop system input-to-state
stability using the descriptor method for Lyapunov-Krasovskii functionals. Due to the trade-off
between delay and gain margin, a limitation of high-gain results for time-delay systems. The limitation
is overcome by using a cascade of sub-predictors. A comparative analysis is also presented, showing that
our proposal achieves a better trade-off between delay and gain margin.
\end{abstract}

\begin{document}

\flushbottom
\maketitle
\thispagestyle{empty}

\section{Introduction}


In delay-systems studies, the control of input delays has received noteworthy attention. The main reason is certainly the widespread 
presence of delays in the control loop of many processes and systems. 
For example, in \cite{hollot2002}, the queue length and capacity and the packages propagation time in transmission control protocol (TCP) routers for active queue management are modeled as an input delay; so is the transmission of information in teleoperated systems, as shown by \cite{nuno2009}. In population studies, for example the observer-predictor for a susceptible, infectious, and recovered (SIR) epidemic model, delays in the application of social distancing measures and in the data recollection process make the model more realistic \citep{castanos2021}. Similarly, a delay is observed in the glucose-insulin model of \cite{borri2017} as the plasma glucose concentration does not change instantly to a change in insulin levels.
While small input delays may often be neglected, larger ones are detrimental to the closed-loop system performance and can cause oscillations or instability.
The interest of the delay community in addressing this challenge is fostered by the fact that the problem is well-suited for the development of specialized efficient control tools. 

A natural idea for mitigating the effect of delays in the control loop is to feed the controller with a prediction or reconstruction of what the state will be at a time instant ahead, 
resulting in the compensation of the delay. A noteworthy advantage of this approach is that existing control laws and corresponding tuning methods for delay-free systems can be employed. 
This idea appeared first in the celebrated frequency domain controller by \cite{smith1957}, whose limitations and solutions in the cases of unstable or uncertain processes are widely reported. 
It is also the core principle of the \emph{finite spectrum assignment} approach \citep{manitus1979}, which relies on a future state reconstruction based on the integral Cauchy formula. 
This method owes its name to the fact that, if the linear system model and the delay are known, and the computation of the integral is exact, the closed-loop spectrum is finite.  However,
the infinite-dimensional nature of the ODE -- delayed-integral-equation interconnection \citep{krstic} resurfaces when uncertainties are present, an issue that is mitigated by filtering the
control law \citep{mondie2003,kharitonov}. A strength of this strategy is that the design phase reduces to choosing a Hurwitz filter matrix.

Another approach, inspired by Luenberger observers, consists in the design of an observer for the state a delay interval ahead of the present time instant \citep{germani2002}. 
The dynamics of the observer consists of a copy of the system model, but the present control law is fed in place of the delayed one. A correction factor depending on the error at the present
time instant is also added \citep{zhou2017}. Attractive features of this method are that the observer can be readily built from the system model, so the approach covers a wide range of classes of systems, such
as nonlinear ones \citep{estrada2017,castanos2021}, and that it also covers the case of systems with partial-state measurements. A disadvantage is that the closed-loop dynamics consists in a
delay system that requires the careful choice of the error gain. Even worse, for large input delays, no gain may exist that ensure closed-loop stability. A remedy to this issue is the introduction
of concatenated sub-predictors that estimate the state a fraction of the delay ahead in time: in this case, both the gain and the number of needed sub-predictors must be found. This task is often
carried out by testing LMI-type conditions for an increasing number of predictors until a solution is found \citep{najafi2013}.

For nonlinear uniformly observable systems, there exists a canonical form written as the sum of two components: a chain of integrators and a nonlinear triangular vector field. In the
delay-free case, a popular scheme for this class of systems is the high-gain observer, a special case of the Luenberger observer in which the error gain can be scaled by tuning a scalar parameter
(the so-called \emph{high gain}). When the scalar parameter is large enough, the nonlinearity entering the error dynamics is dominated, thus ensuring a stable asymptotic closed-loop 
behavior \citep{gauthier1992,gauthier1994,gauthier}. The extension to nonlinear time-delay systems was studied by \cite{farza2010, ahmed2012,lei2016a}. It is observed that, in contrast 
with the delay-free case, the scalar gain is limited by the delay. This trade-off imposes restrictions either on the size of the delay or the size of the Lipschitz constant characterizing
the nonlinearities. The typical design approach consists in tuning the observer for the delay-free system and then treating the delay as an additional perturbation.

The aim of this contribution is to reduce conservatism in the above-mentioned trade-off. First, high-gain predictors are analyzed under the lens of homogeneity \citep{bacciotti}. 
Instead of treating the delay as a perturbation, the observation error dynamics are directly tuned as a delay system using the \emph{multiplicity-induced dominance} (MID) tools recently developed
\citep{balogh2022}, especially those addressing this particular case \citep{rojas2022b}. The delay margin of the linearized error dynamics is computed exactly using frequency
domain analysis, while the gain margin is computed using an \emph{input-to-state stability} (ISS) perspective. Finally, the homogeneity point of view permits us integrate the delay and
gain margins into a less conservative trade-off.

The paper is organized as follows: Section~\ref{sec:problem} gives a brief overview of the class of systems under consideration, of the prediction problem, and relevant assumptions. Technical facts on weighted homogeneity and input-to-state stability are also introduced.
In Section~\ref{sec:predictor}, the high-gain prediction approach is presented. The MID tuning methodology for the predictor gain ensuring the convergence of the prediction error dynamics is summarized. The trade-off between delay and gain margin arising when applying the high-gain approach to delay systems is discussed. An extension to a cascade of sub-predictors is shown to overcome the limitations produced by this trade-off. 
In Section~\ref{sec:discussion}, we review the high-gain predictors introduced by \cite{ahmed2012} and \cite{lei2016a}, and compare them with our proposal, which appears to achieve a less restrictive trade-off.
We end the contribution with an academic example in Section~\ref{sec:example}, and finally provide some concluding remarks.

\section{Preliminaries and Problem Formulation} \label{sec:problem}

In this section, we state the problem formally and recall a few basic notions on homogeneity and 
input-to-state stability in the framework of time-delay systems.

\subsection{Notation}

Here, $\mathbb{N}$, $\RE$, and $\REp = \left\{ \theta \in \RE \mid \theta \ge 0 \right\}$ stand
for the sets of natural, real, and nonnegative real numbers, respectively.
The space of $\RE^n$-valued continuous functions of a real variable is denoted by $\CO(\RE,\RE^n)$. Likewise,
$\LE(\RE,\RE^m)$ denotes the space of $\RE^m$-valued Lesbegue-measurable functions of a real variable.
We adopt the standard notation for systems with time delays: 
Given $x \in \CO(\RE,\RE^n)$ and $h \ge 0$, we use $x_t \in \COh$ to denote the function defined 
by $x_t(\theta) = x(t+\theta)$ for $-h \le \theta \le 0$. The function $u_t \in \LE([-h,0],\RE^m)$ is defined 
similarly for $u \in \LE(\RE,\RE^m)$. Given $\xi \in \RE^n$, $\phi \in \COh$, and $u \in \LE(\RE,\RE^m)$,
we denote by
\begin{displaymath}
	\norm{\xi} \;, \quad \norm{\phi}_h = \max_{\theta \in [-h,0]} \norm{\phi(\theta)} \;, \quad \text{and} \quad
		\norm{u}_\infty = \esssup_{t \in \RE} \norm{u(t)} \;,
\end{displaymath}
respectively, the Euclidean norm, the standard uniform norm, and the $\LE^\infty$ norm.


Given $k,n \in \mathbb{N}$ such that $k \le n$, `$n$ choose $k$' is indicated by
\begin{displaymath}
    \binom{n}{k} = \frac{n!}{k!(n-k)!} \;.
\end{displaymath}

A function $\alpha \in \CO(\REp,\REp)$ is said to be of class $\CK$ if it is strictly
increasing and $\alpha(0) = 0$. It is said to be of class $\CK_\infty$ if $\alpha \in \CK$ and
$\lim_{\theta\to+\infty} \alpha(\theta) = +\infty$. A function 
$\beta : \REp\times \REp \to \REp$ is said to be of class $\CKL$ if, for 
each fixed $t \ge 0$, $\beta(\cdot,t) \in \CK$ and, for each fixed $\theta \in \REp$, $\beta(\theta,\cdot)$
is continuous, non-increasing, and $\lim_{t \to +\infty}\beta(\theta,t) = 0$.

\subsection{Problem statement}

Consider an open, connected subset $\Omega \subset \RE^n$ and a control-affine system with delayed inputs of the form
\begin{subequations} \label{eq:normal_form}
\begin{align}
        \dot{x}(t) &= Ax(t) + \varphi(x(t),u(t-h)) \label{eq:system} \\
              y(t) &= Cx(t) \;, \label{eq:output}
\end{align}
where $x(t) \in \Omega$, $u(t) \in U \subset \RE^m$, and $y(t) \in \RE$ are the state, the input, and the output
at time $t \in \REp$. The nonnegative constant $h$ is the input delay. The matrix $A \in \RE^{n\times n}$
consists of a single Jordan block with eigenvalue 0, $C = [1,0,\ldots,0] \in \RE^{1\times n}$, and 
$\varphi : \Omega\times U \to \RE^n$ is triangular in $x$,
\begin{equation} \label{eq:phi}
    \varphi(x,u) = 
    \begin{bmatrix}
        \varphi_1(x_1,u) & \varphi_2(x_1,x_2,u) & \cdots & \varphi_n(x_1,\ldots,x_n,u)
    \end{bmatrix}^\top \;,
\end{equation}
\end{subequations}
and affine in $u$. The domain $\Omega$ is supposed to be invariant with respect to~\eqref{eq:normal_form}.

Recall that an input is said to be \emph{universal} if the resulting outputs are not identical
for any pair of different initial conditions. A system is \emph{uniformly observable} if every
input is universal (see~\cite{besancon} for details). It turns out that systems of the 
form~\eqref{eq:normal_form} are uniformly observable. Moreover, by a change of coordinates, 
any uniformly observable system can be transformed into the canonical form~\eqref{eq:normal_form},
provided it is sufficiently smooth \citep{gauthier1981}. The latter property makes the 
class~\eqref{eq:normal_form} a reasonable target for observation problems.

Following \cite{gauthier1981}, we further impose the following.

\begin{ass} \label{ass:Lipschitz}
    The vector field~\eqref{eq:phi} can be extended to the whole $\RE^n$ space by a function that is globally Lipschitz in $x$, uniformly in
    $u$. In other words, there exists a nonnegative constant $\gamma_\varphi$ such that
    \begin{equation*}
        \norm{\varphi(\hat{x},u) - \varphi(x,u)} \leq \gamma_\varphi \norm{\hat{x}-x}
    \end{equation*}
	for all $\hat{x},x \in \RE^n$ and all $u \in U$.
\end{ass}

Global Lipschitz continuity in $\Omega$ is a standard assumption ensuring the uniqueness of solutions, while the
possibility to extend the vector field to $\RE^n$ without destroying globality
is a rather benign addition. The assumption admits the use of linear terms to dominate the nonlinearities.

\begin{prob} \label{prob:prediction}
    Use $u_t \in \LE([-h,0],\RE^m)$ and $y_t \in \CO([-h,0],\RE)$ to find an estimate
    $\hat{x}(t)$ for $x(t + h)$ such that
    \begin{equation*} \label{eq:prediction}
        \lim_{t \to +\infty}\left( \hat{x}(t)-x(t+h) \right) = 0.
    \end{equation*}
\end{prob}

The motivation is clear: Suppose there exists a feedback law $\hat{u}:\RE^n \to \RE^m$ such that
the target system
\begin{equation} \label{eq:target}
	\dot{x}(t) = Ax(t) + \varphi(x(t),\hat{u}(x(t))
\end{equation}
satisfies the design specifications (e.g., asymptotic stability of the origin). Because we have only partial state
information and because the input is retarded, such a system is not feasible. However, it is reasonable to expect the
system $\dot{x}(t) = Ax(t) + \varphi(x(t),\hat{u}(\hat{x}(t-h))$ --- which is indeed attainable --- to behave `closely' 
to~\eqref{eq:target}. There are different ways to formalize the previous statement (e.g., as carried out by
\cite{lei2016a, estrada2017}), but we will not do so in this paper. Instead, we will focus solely on the prediction problem.

In principle, vector outputs $y(t) \in \RE^o$ can also be considered in the problem statement \citep{hammouri}.
Indeed, when $o > m$, the class of uniformly observable systems becomes more interesting, as uniform observability
becomes a generic property \citep{gauthier}. Also, using the Fundamental Theorem of Calculus, it is
possible to incorporate time-varying delays $h + \delta h(t)$. The usual approach is to write the delayed error as
$e(t-h-\delta h(t)) = e(t-h) - \int_{t-h-\delta h(t)}^{t-h} \dot{e}(s)\rd s$ and treat the integral term
as a perturbation. However, the notation becomes cumbersome and distracts our attention from the main results, so we 
prefer to avoid the multi-output and time-varying delay scenarios.

We will approach Problem~\ref{prob:prediction} within the frameworks of homogeneity and input-to-state
stability, two subjects that we now briefly revisit while including a couple of new results. 
The fact that~\eqref{eq:normal_form} has a retarded input
invites us to consider the time-delayed versions of such frameworks.

\subsection{A homogeneity perspective of high-gain observers}

Allow us to recall the notion of weighted homogeneity (see \cite{bacciotti} for a detailed exposition).

\begin{defn}
    Let $r = [r_1, \ldots, r_n]$ be an $n$-tuple of positive real numbers.
    \begin{itemize}
        \item The one-parameter family of \emph{dilations} $\Lambda^r_{\lambda}:\RE^n \to \RE^n$ is defined by
         \begin{equation*}
             \Lambda^r_\lambda(x) = 
             \begin{bmatrix}
                \lambda^{r_1}x_1 & \cdots & \lambda^{r_n}x_n 
             \end{bmatrix}^\top
             \;, \quad x \in \RE^n \;, \quad \lambda > 0 \;.
         \end{equation*}
         The numbers $r_i$ are the \emph{weights}.
        \item A function $V : \RE^n \to \RE$ is said to be $r$-\emph{homogeneous} of degree $d \in \RE$ if 
         \begin{equation*}
             V\left(\Lambda^r_\lambda(x)\right) = \lambda^d \cdot V(x) \;, \quad x \in \RE^n \;, \quad \lambda > 0 \;.
         \end{equation*}
        \item A vector field $f : \RE^n \to \RE^n$ is said to be $r$-\emph{homogeneous} of degree $d \in \RE$ if 
         \begin{equation} \label{eq:homo_vector}
             f\left(\Lambda^r_\lambda(x)\right) = \lambda^d \cdot \Lambda^r_\lambda\left(f(x)\right) \;, 
					\quad x \in \RE^n \;, \quad \lambda > 0 \;.
         \end{equation}       
    \end{itemize}
\end{defn}

Note that the composition of dilations gives $\Lambda_\mu^r\circ\Lambda_\lambda^r = \Lambda_{\mu\cdot\lambda}^r$ for all 
$\lambda, \mu > 0$, so the family of dilations constitutes a one-parameter commutative group of linear transformations
with identity $\Lambda_1^r$ and inverse $(\Lambda^r_\lambda)^{-1} = \Lambda_{\lambda^{-1}}^r$.

\begin{lem} \label{lem:homo_linear}
    Let $A$ and $C$ be as in~\eqref{eq:normal_form}. The vector field $A x$ and the function $C x$ are
    $r$-homogeneous of degree $d=1$ with $r = [1,\ldots,n]$.
\end{lem}

One way to interpret~\eqref{eq:homo_vector} is that, modulo a time rescaling, $f$ is invariant
to coordinate transformations induced by dilations.
\begin{prop} \label{prop:dilation}
    Let $t \mapsto x(t)$ be an integral curve of $f : \RE^n \to \RE^n$ with $f$ an $r$-homogeneous vector field of degree $d$.
    Then, for any $\lambda > 0$,
\begin{equation*}
    \xi : t \mapsto \Lambda_\lambda^r\left(x(\lambda^{d}\cdot t)\right)
\end{equation*}
is also an integral curve of $f$.
\end{prop}

An interesting consequence of Proposition~\ref{prop:dilation} is that local qualitative properties of homogeneous vector
fields are inherently global. In particular, asymptotic stability implies global asymptotic stability~\citep{bacciotti}.

While the majority of systems are not naturally defined by homogeneous vector fields, many can be represented as a formal series 
of homogeneous ones~\citep{hermes1991}. It is also possible to approximate non-homogeneous systems by
homogeneous vector fields~\citep{bacciotti}. The most common approximations are, of course, linear vector fields, which are all
$r$-homogeneous of degree $d=0$ with the \emph{standard weights} $r=[1,\ldots,1]$. However, approximation by systems of
degrees different from zero has also proved useful for, e.g., performing nonlinear controllabilty
tests~\citep{hermes1991}.



The following shows that while $\varphi$ is not necessarily homogeneous, its Lipschitz continuity is not destroyed by dilation transformations.

\begin{lem} \label{lem:Lip_homogeneous}
    Consider the controlled vector field~\eqref{eq:phi} and define 
    \begin{equation*}
        \Delta\varphi(e,x,u) = \varphi(x+e,u) - \varphi(x,u)
    \end{equation*}
    and 
    \begin{displaymath}
        \Delta\Phi_\lambda\left(\varepsilon,x,u\right) = 
         \Lambda_{\lambda^{-1}}^{r}\left( \Delta\varphi\left(\Lambda_\lambda^{r}(\varepsilon),x,u\right) \right) \;.
    \end{displaymath}
    For $r = [1,\ldots,n]$, Assumption~\ref{ass:Lipschitz} implies that there exists a nonnegative constant
    $\gamma_{\Phi}$ such that
    \begin{equation} \label{eq:lip_l}
        \norm{\Delta\Phi_\lambda\left(\varepsilon,x,u\right)} \le \gamma_{\Phi} \norm{\varepsilon}
    \end{equation}
    for all $x \in \Omega$, $\varepsilon \in \RE^n$, $u \in U$, and $\lambda \ge 1$.
\end{lem}


Homogeneity can be used as a design tool. For example, in the usual approach to higher-order sliding-mode control, the plant is
rendered homogeneous by means of a homogeneous feedback~\citep{levant2005}. Also, dilations can be explicitly 
incorporated in the feedback control law so that the feedback gain is characterized by a single parameter $\lambda$. This 
strategy has been used to construct time-varying feedback gains that achieve improved convergence~\citep{mcloskey1997} or 
even prescribed-time stability~\citep{chitour2020}, and to construct high-gain observers with adaptive gains and homogeneous 
correction terms~\citep{andrieu2009}.

The notion of homogeneity has been extended to functionals and functional vector fields. Many --- albeit not all --- properties 
of delay-free homogeneous systems can be translated to homogeneous systems with time delays. It has been shown, for example, that
a time-delay version of Proposition~\ref{prop:dilation} holds \citep{efimov2014}.

The tool of choice for estimating $x(t)$ is the high-gain observer \citep{hammouri}
\begin{equation} \label{eq:high_gain_obs}
    \dot{\hat{x}}(t) = A\hat{x}(t) + \Lambda_\lambda^r\left(L\cdot\left(y(t)-C\hat{x}(t)\right)\right) + 
     \varphi(\hat{x}(t),u(t)) \;,
\end{equation}
where $L = [l_1,\ldots,l_n]^\top \in \RE^n$ and $\lambda > 0$ are design parameters. The well-known advantages 
of~\eqref{eq:high_gain_obs} are summarized as a corollary to Proposition~\ref{prop:high_gain_pred} 
(given below in Section~\ref{sec:predictor}).

\begin{cor} \label{cor:high_gain_obs}
    Consider a uniformly observable system~\eqref{eq:normal_form} and a high-gain observer~\eqref{eq:high_gain_obs}. 
    Define the time-scale $\tau = \lambda t$ and the \emph{dilated error}
    \begin{equation} \label{eq:dilated_error}
      \varepsilon(\tau) = \Lambda_{\lambda^{-1}}^r\left(e(\lambda^{-1}\tau)\right)  
    \end{equation}
    with $e(t) = x(t) - \hat{x}(t)$ and $r = [1,\ldots,n]$. Let 
    $\Upsilon_{\lambda}^L(e) = Ae - \Lambda_{\lambda}^r\left(LC e\right)$.
    The dilated error evolves according to
    \begin{equation} \label{eq:finite_error}
        \frac{\rd\varepsilon(\tau)}{\rd \tau} = \Upsilon_{1}^L\left(\varepsilon(\tau)\right) 
        + \lambda^{-1}\cdot\Delta\Phi_\lambda\left(\varepsilon(\tau),\hat{x}(\lambda^{-1}\tau),u(\lambda^{-1}\tau-h)\right) \;.
    \end{equation}
\end{cor}

\begin{rem} \label{rem:observer}
    Equation~\eqref{eq:finite_error} suggests that the observation problem be solved in two steps:
    \begin{enumerate}
        \item Choose $L$ so that $\Upsilon_{1}^L(\varepsilon) = (A-LC)\varepsilon$ is exponentially stable. This is always
         possible by the observability of the pair $(A,C)$.
        \item Simply choose $\lambda$ large enough so that $\lambda^{-1}\cdot \Delta\Phi$ is dominated by $\Upsilon_{1}^L$.
         This is always possible by Lemma~\ref{lem:Lip_homogeneous}.
    \end{enumerate}
\end{rem}

Our primary strategy is to transfer Remark~\ref{rem:observer} (about the standard observation problem) to 
the more complicated Problem~\ref{prob:prediction}.

\subsection{Input-to-state stability for systems with time delays}

We also consider general time-delay systems of the form
\begin{equation} \label{eq:general}
    \dot{e}(t) = f(e_t,w(t)) \;,
\end{equation}
where $f : \CO([-h,0],\RE^{n_1})\times\RE^p \to \RE^{n_1}$ is a `functional vector field'. The function
$e_t \in \CO([-h,0],\RE^{n_1})$ and the vector $w(t) \in \RE^p$ are, respectively, the state and input
at time $t$. We assume that the reader is familiar with basic properties such as the existence and uniqueness of solutions,
stability, and convergence in the time-delay setting \citep{kharitonov}.

Input-to-state stability \citep{sontag} extends the notion of Lyapunov stability to systems with inputs. 
During the past decades, the concept has been developed for systems with time delays, and it now provides a suitable
framework for solving Problem~\ref{prob:prediction}. We briefly recall the basic notions. For a broader view,
we recommend the recent survey by \cite{chaillet2023}.

\begin{defn}
    System~\eqref{eq:general} is said to be \emph{input-to-state stable (ISS)} if there exists functions
    $\beta \in \CKL$ and $\mu \in \CK$ such that, for all $e_0 \in \CO([-h,0],\RE^{n_1})$ and all
    $w \in \LE(\REp,\RE^p)$,
    \begin{displaymath}
        \norm{e(t)} \le \beta(\norm{e_0}_h, t) + \mu(\norm{w}_\infty) \;, \quad t \in \REp \;.
    \end{displaymath}
    The function $\mu$ is called an \emph{ISS gain}.
\end{defn}

ISS implies the global asymptotic stability of the input-free system and the boundedness of solutions in the face of 
bounded inputs. In analogy with the original delay-free setting, linear time-delay systems that are asymptotically stable are
inevitably ISS.

\begin{cor}[to Prop. 2.5 in \cite{pepe2006}] \label{cor:linear}
    Consider the functional $f(\phi,w) = \psi(\phi) + Bw$ with $\psi$ a linear bounded operator from
    $\CO([-h,0],\RE^{n_1})$ to $\RE^{n_1}$ and $B \in \RE^{{n_1}\times p}$. If the origin of
    $\dot{e}(t) = f(e_t,0)$ is asymptotically stable, then~\eqref{eq:general} is ISS.
    Moreover, the ISS gain can be chosen to be linear.
\end{cor}

One of the benefits of the ISS framework is the extensive repertoire of 
Lyapunov characterizations.

\begin{defn} 
    A functional $V:\COh \to \RE$ is said to be \emph{Lipschitz on bounded sets}
    if, for any $\eta \in \REp$, there exists $\gamma > 0$ such that
    \begin{displaymath}
        |V(\phi_1) - V(\phi_2)| \le \gamma \norm{\phi_1-\phi_2}_h 
    \end{displaymath}
    for all $\phi_1, \phi_2 \in \COh$ with $\norm{\phi_1}_h \le \eta$ and $\norm{\phi_2}_h \le \eta$.
\end{defn}

While taking the derivative of a functional $V:\COh \to \RE$ along a system trajectory, we will
make use of \emph{Driver's derivative} $D^+V : \COh\times\RE^n \to \RE$, defined by
\begin{displaymath}
    D^+ V(\phi,v) =  \limsup_{\Delta\theta \to 0^+} \frac{V(\phi_{\Delta\theta,v})-V(\phi)}{\Delta\theta} \;,
\end{displaymath}
where
\begin{displaymath}
    \phi_{\Delta\theta,v} =
    \begin{cases}
     \phi(\theta+\Delta\theta),             & \text{if $\theta \in [-h,-\Delta\theta)$} \\
     \phi(0) + v\cdot(\theta+\Delta\theta), & \text{if $\theta \in [-\Delta\theta,0]$}
    \end{cases}
\end{displaymath}
(see \cite{chaillet2023} for details).

\begin{defn}
    A functional $V:\COh \to \REp$ is said to be a \emph{Lyapunov-Krasovskii functional candidate (LKF)} if it is
    Lipschitz on bounded sets and there exists $\alpha_1, \alpha_2 \in \CK_\infty$ such that, for all
    $\phi \in \COh$,
    \begin{displaymath}
        \alpha_1(\norm{\phi(0)}) \le V(\phi) \le \alpha_2(\norm{\phi}_h) \;.
    \end{displaymath} 
\end{defn}

\begin{thm}[Thm. 2 in \cite{kankanamalage2017}] \label{thm:ISS}
    System~\eqref{eq:general} is ISS if, and only if, there exists a LKF $V$ and functions $\alpha_3, \chi \in \CK$
    such that 
    \begin{equation} \label{eq:LKF_ISS}
        V(\phi) \ge \chi(\norm{v}) \quad \Longrightarrow \quad D^+ V(\phi,f(\phi,v)) \le -\alpha_3(\norm{\phi(0)}) \;. 
    \end{equation}
    for all $\phi \in \COh$.
\end{thm}

ISS has been found to be a useful concept to study, e.g., robustness and stability properties
of nonlinear interconnected systems. For example, ISS implies robustness with respect to bounded delay-free feedback.

\begin{cor} \label{cor:robust}
    Let~\eqref{eq:general} be ISS with input $w$. There exists a function $\rho \in \CK_\infty$, called a 
    \emph{gain margin} such that, for every feedback law $\kappa : \RE^{n_1}\times\RE \to \RE^p$ satisfying
    $\norm{\kappa(\xi,t)} \le \rho\left( \norm{\xi} \right)$, the system
    \begin{equation} \label{eq:feed_pert}
        \dot{e}(t) = f(e_t,\kappa(e_t(0),t)+w(t))
    \end{equation}
    is again ISS with input $w$.
\end{cor}

For linear systems, gain margins can be estimated by solving linear matrix inequalities.

\begin{prop} \label{prop:margin}
    Consider a linear time-delay system of the form $\dot{e}(t) = Ae_t(0) + A_1e_t(-h) + w(t)$
    with $A, A_1 \in \RE^{n\times n}$ and $w(t) \in \RE^n$. The system is ISS with a gain margin 
    $\rho(\theta) = \gamma_m \cdot \theta$ if there exist matrices $P>0$, $R>0$, $S>0$, $P_2$, $P_3$, $P_4 \in \RE^{n\times n}$
    such that $W < 0$ with
	\begin{equation} \label{eq:W}
		W = 
         \begin{bmatrix}
				A^\top P_2 + P_2^\top A+S-R+\gamma_m^2 I_{n} &  P-P_2^\top+A^\top P_3 & P_2^\top A_1+R &  P_2^\top + A^\top P_4 \\
				                                           * & -P_3-P_3^\top + h^2R   & P_3^\top A_1   &  P_3^\top - P_4 \\
				                                           * & *                      & -S-R            & A_1^\top P_4 \\
				                                           * & *                      & *               &  P_4^\top + P_4 - I_{n}
			\end{bmatrix} \;.
		\end{equation}
    Here, the asterisks correspond to the unique terms that make $W$ symmetric. 
\end{prop}

Since ISS prevents finite escape time, ISS is closed under cascade interconnections. To formulate this formally, 
consider another general system
\begin{equation} \label{eq:general_2}
    \dot{z}(t) = g(z_t,\overline{w}_t) \;,
\end{equation}
where $g : \CO([-h,0],\RE^{n_2})\times\LE([-h,0],\RE^{n_1}) \to \RE^{n_2}$, $z_t \in \CO([-h,0],\RE^{n_2})$ and 
$\overline{w} \in \LE([-h,0],\RE^{n_1})$.

We will make use of the following.

\begin{cor}[\cite{chaillet2023}] \label{cor:cascade}
    Suppose that~\eqref{eq:general} and~\eqref{eq:general_2} are ISS with inputs $w$ and $\overline{w}$, respectively.
    Then, the cascade resulting from the interconnection $\overline{w} = e$ is again ISS with input $w$.
\end{cor}

\section{High-gain sub-predictors} \label{sec:predictor}

In the spirit of~\cite{germani2002}, we solve the problem of predicting future state values
with a Luenberger-like observer, that is, by constructing a copy of the plant and adding output
injection to stabilize the error dynamics. However, unlike a true Luenberger observer, the output injection 
of the predictor is delayed by $h$ units of time. To cope with the nonlinear term $\varphi$, we will use a 
high-gain structure injection similar to~\eqref{eq:high_gain_obs}. 

Let $\bar{h} \in [0,h]$ and consider a high-gain predictor~\citep{ahmed2012} of the form
\begin{equation} \label{eq:high_gain_pred}
    \dot{\hat{x}}(t) = A\hat{x}(t) + \Lambda_\lambda^r\left(L\cdot\left(y(t)-C\hat{x}(t-\bar{h})\right)\right) + 
     \varphi(\hat{x}(t),u(t-h+\bar{h})) \;,
\end{equation}
where $L \in \RE^n$ and $\lambda > 0$ are again design parameters.

It can be readily seen that the error
\begin{equation} \label{eq:error}
    e(t) = \hat{x}(t-\bar{h}) - x(t) 
\end{equation}
evolves according to
\begin{equation} \label{eq:error_der}
    \dot{e}(t) =  \psi_{\lambda,\bar{h}}^L\left(e_t\right) + \Delta\varphi(e(t),x(t),u(t-h)) \;,
\end{equation}
where $\psi_{\lambda,\bar{h}}^L : \CO([-\bar{h},0],\RE^n) \to \RE^n$ is given by
\begin{equation} \label{eq:psi}
    \psi_{\lambda,\bar{h}}^L\left(e_t\right) = Ae_t(0) - \Lambda_\lambda^r\left(LCe_t(-\bar{h})\right) \;.
\end{equation}

\begin{lem} \label{lem:pseudo_homo} 
    Let $A$ and $C$ be as in~\eqref{eq:normal_form}. For any $L \in \RE^n$, the functional vector field~\eqref{eq:psi} 
    satisfies
    \begin{equation} \label{eq:pseudo_homo}
        \psi_{\lambda,\bar{h}}^L\left( \Lambda_\lambda^r\circ e_t \right) = 
          \lambda\cdot\Lambda_{\lambda}^r\left(\psi_{1, \bar{h}}^L(e_t)\right) \;, \quad 
            e_t \in \CO([-\bar{h},0],\RE^n) \;, \quad \lambda > 0 \;.
    \end{equation}
\end{lem}

The property described in Lemma~\ref{lem:pseudo_homo} (cf. the homogeneity property~\eqref{eq:homo_vector}) is key to the following.

\begin{prop} \label{prop:high_gain_pred}
    Consider a uniformly observable system~\eqref{eq:normal_form}, a high-gain predictor~\eqref{eq:high_gain_pred}, and
    the dilated error~\eqref{eq:dilated_error} with $e(t)$ as in~\eqref{eq:error}, $r = [1,\ldots,n]$, and $\tau = \lambda t$.
    The dilated error evolves according 
    to the time-delay differential equation
    \begin{equation} \label{eq:finite_error_delay}
        \frac{\rd\varepsilon(\tau)}{\rd \tau} = \psi_{1,\lambda \bar{h}}^L\left(\varepsilon_\tau\right) 
        + \lambda^{-1}\cdot\Delta\Phi_\lambda\left(\varepsilon_\tau(0),x_{\lambda^{-1}\tau}(0),u_{\lambda^{-1}\tau}(-h)\right) \;.
    \end{equation}
\end{prop}

We will henceforth set $\bar{h} = h$ for simplicity.

\begin{rem} \label{rem:predictor}
Just as with Remark~\ref{rem:observer}, $\lambda^{-1}\cdot\Delta\Phi_\lambda$ in~\eqref{eq:finite_error_delay} can be dominated by 
setting $\lambda$ large enough. Unfortunately, when $h > 0$, the functional vector field $\psi_{1,\lambda h}^L$ is no 
longer independent of $\lambda$, so the problem cannot be easily solved in two separate steps as before. First,
it is quite possible that $\psi_{1,\lambda h}^L$ is exponentially stable for some values of
$\lambda$ and unstable for others. Second, the spectrum of $\psi_{1,\lambda h}^L$ is infinite, so designing
$L$ such that $\psi_{1,\lambda h}^L$ is exponentially stable is no longer a trivial task, even if $\lambda$ is fixed.
\end{rem}

The usual approach to solving the problem outlined in Remark~\ref{rem:predictor} is first choosing $L$ so that $A-LC$ is Hurwitz,
effectively considering the system as delay-free. In a subsequent stage, the true delay is treated as a perturbation, and a set
of conditions are imposed on how small should $h$ and $\gamma_\Phi$ be to preserve asymptotic stability (see, e.g., 
\cite{ahmed2012,najafi2013,lei2016a}). 
Needless to say, considering both $\varphi$ and $h$ as perturbations results in
rather conservative conditions (see Section~\ref{sec:discussion} for more details). In what follows, we present a 
methodology for designing $L$ so that the spectrum of $\psi_{1,\lambda h}^L$ is located in the open left-half plane,
regardless of how big $h$ is. Then, tuning $\lambda$ in a less conservative fashion is relatively straightforward.

\subsection{Designing the vector gain $L$} \label{sec:L}

We focus on the functional vector field $\psi_{1,\lambda h}^L$ and recall some results from~\cite{rojas2022b}.
It is straightforward to verify that the characteristic quasipolynomial 
$D_{\lambda h}^L(s) = \det\left( sI - A - \eu^{-\lambda h s}LC \right)$ associated to $\psi_{1,\lambda h}^L$ is
\begin{equation} \label{eq:D}
    D_{\delta}^L(s) = s^n + \left( l_1 s^{n-1} + l_2 s^{n-2} + \cdots + l_n \right) \eu^{-\lambda h s} \;,
\end{equation}
where $\delta = \lambda h$ is nothing but the delay in the $\tau$ time scale.

\begin{defn}
    We say that a root $s^\star \in \mathbb{C}$ of $D_{\delta}^L(s)$ is \emph{dominant} if
    \begin{displaymath}
        \Re (s^\star) = \max\left\{ \Re(s) \mid D_{\delta}^L(s) = 0 \right\} \;.
    \end{displaymath}
\end{defn}

We begin by characterizing all the roots of $D_{\delta}^L(s)$ that can be multiply assigned by $L$.  Finding a root with multiplicity $n+1$ is equivalent \citep{michiels} to find a solution $s$ for
\begin{equation} \label{eq:derivadasDs}
    D_\delta^L(s)\Big|_{s=\sigma} = 0 \;, \quad \frac{\rd D_\delta^L(s)}{\rd s}\Big|_{s=\sigma} = 0 \;, \quad \cdots,
  \quad \frac{\rd^n D_\delta^L(s)}{\rd s^n}\Big|_{s=\sigma} = 0 \;.
\end{equation}
In our special case, the notation introduced by \cite{balogh2022} allows an elegant rewriting of conditions~\eqref{eq:derivadasDs} in matrix form. For  each $k=0,1,\ldots,n$, we implicitly define $R_k(s, \delta)$,
a polynomial in $s$ and $\delta$ given by
\begin{equation} \label{eq:Rk_def}
    R_k(s,\delta)e^{s \delta} = \frac{\rd^k P(s)e^{s \delta}}{\rd s^k} \;.
\end{equation}
Multiplying~\eqref{eq:derivadasDs} by $e^{s}$  and using the implicit definition~\eqref{eq:Rk_def}, 
we can express~\eqref{eq:derivadasDs} as
\begin{subequations} \label{eq:Matrix_form}
\begin{align}
 M(s) J_n L + R(s,\delta)e^{s} &= 0 \label{eq:Matrix_form_derivatives} \\
	    R_n(s,\delta)e^{s} &= 0 \label{eq:Matrix_form_polynomial}
\end{align}
\end{subequations}
with 
\begin{displaymath}
 R(s,\delta) = 
 \begin{bmatrix}
  R_0(s,\delta) \\ R_1(s,\delta) \\ \vdots \\ R_{n-1}(s,\delta)
 \end{bmatrix}, \quad J_n = 
 \begin{bmatrix}
  0&\cdots&0&1 \\ 
  0&\cdots&1&0 \\ 
  \vdots& &\vdots&\vdots \\ 
  1&\cdots&0&0
 \end{bmatrix},
\end{displaymath}
where $M(s) = \left[m_{ij}(s)\right]$ is a $n\times n$ matrix whith polynomials entries in $s$
\begin{displaymath}
 m_{ij}(s) = 
 \begin{cases}
  \frac{(j-1)!}{(j-i)!}s^{j-i} & \text{if $j \ge i$} \\
  0                            & \text{if $j < i$}
 \end{cases} \;, \quad i,j = 1,\ldots,n \;
\end{displaymath}
and $R_k(s,\delta)$ takes the explicit form
\begin{equation} \label{eq:Rk_q}
 R_k(s,\delta) := \sum_{i=1}^{k+1}\binom{n}{k-i+1} \frac{k!}{(i-1)!}s^{n-k+i-1}\delta^{i-1}. 
\end{equation}

It turns out that $R_n(s,\delta)$ in \eqref{eq:Matrix_form_polynomial} is independent of $L$ and reduces to a polynomial in $\sigma=\delta s$. Thus, it has $n$ easy-to-compute roots corresponding to roots of multiplicity $n+1$ of \eqref{eq:D}. By Sturm's Theorem, these  $n$ roots are shown to be real and negative. 

\begin{lem}[\cite{rojas2022b}] 
    The polynomial
    \begin{equation} \label{eq:q}
        q(\sigma) := \sum_{j=0}^{n} \binom{n}{j} \frac{n!}{j!} \sigma^{j} \;.
    \end{equation}
    has $n$ negative and distinct real roots.
\end{lem}

In contrast with the work of \cite{balogh2022}, where an $n^{\mathrm{th}}$-degree quasipolynomial with $2n$ free parameters lead to a single solution which is evidently dominant, here we must determine which of the solutions of $q(\sigma) = 0$ corresponds to the dominant root of multiplicity $n+1$ of (\ref{eq:D}). The careful analysis of the integral term in the factorization inspired by \cite{mazanti2021}, which holds when $D_{\delta}^{L}(s)$ has a root of multiplicity $n+1$,
\begin{equation} \label{eq:Ds_n1_roots}
  D_\delta^L(s) = \frac{1}{n!} \left(s-\sigma\right)^{n+1} \int_{0}^{\delta}R_n(\sigma,\theta)e^{-(s-\sigma)\theta}\rd\theta \;, 
\end{equation}
shows that it is the rightmost root of the polynomial (\ref{eq:q}). Finally, the explicit gains assigning the dominant root of multiplicity $n+1$ are obtained by substituting the rightmost root of \eqref{eq:q} into \eqref{eq:Matrix_form_derivatives} and solving for $L$.

\begin{thm}[\cite{rojas2022b}] \label{thm:L} 
    Let $L^\star = [l_1^\star,\ldots,l_n^\star]^\top \in \RE^n$ be defined by
    \begin{equation*}
        l^\star_k = \sum_{j=n-k+1}^{n}\sum_{i=1}^{j} \left[(-1)^{n+j+k}\binom{n}{j-i} \binom{j-1}{n-k} 
         \frac{1}{(i-1)!}{\sigma_\star}^{i-1}\right]{\sigma_\star}^k \eu^{{\sigma_\star}} \;, \quad k=1,\ldots,n \;,
    \end{equation*}
    where $\sigma_\star$ is the rightmost root of polynomial~\eqref{eq:q}. For $L = \Lambda_{1/\delta}^r(L^\star)$, $D_{\delta}^{L}(s)$ has a dominant root of multiplicity $n+1$ at $\sigma_\star/\delta$.
\end{thm}

Notice that, for $L = \Lambda_{1/\delta}^r(L^\star)$, we have $\psi_{1,\delta}^L = \psi_{1/\delta,\delta}^{L^\star}$
and that, by Lemma~\ref{lem:pseudo_homo}, the functional vector field is equivalent (modulo a time and state-space re-scaling) to
$\psi_{1,1}^{L^*}$. Thus, for the purpose of assigning a dominant multiple root, there is no loss of generality in the restriction
$\delta = 1$ and the corresponding use of the gain $L=L^\star$.

It follows from the previous discussion that the functional vector field $\psi_{1,\lambda h}^{L^\star}$ with $\lambda = 1/h$ is stable with an 
exponential convergence $\sigma_\star$ (in the original time-scale $t$, the convergence rate is
$\lambda \sigma_\star = \sigma_\star/h$). However, it is not clear at this point whether or not $\lambda = 1/h$ is large enough
for $\psi_{1,1}^{L^\star}$ to dominate $\lambda^{-1}\cdot \Delta\Psi_\lambda$ in~\eqref{eq:finite_error_delay}.
We address this question in the following subsections.

\subsection{Stability intervals for the scalar gain $\lambda$} \label{sec:delta}

Now we fix $L = L^\star$ and ask ourselves: `For which values of $\lambda$ (other than $1/h$) is
$\psi_{1,\lambda h}^{L^\star}$ exponentially stable?' Recall that $\delta = \lambda h$ is nothing but the delay
in the $\tau$ time scale. Thus, to answer the question above, we will analyze the stability regions of $\psi_{1,\delta}^{L^\star}$
in the delay parameter space. Afterward, any stable or unstable value for $\delta$ can be mapped back, respectively,
to a stable or unstable value $\lambda=\delta/h$. 

To determine the stability regions, we will use a $\delta$-decomposition method ($\tau$-decomposition method in the book 
by \cite{michiels}).
In accordance with the method, we will partition the delay space into open intervals,
\begin{equation} \label{eq:partition}
    \REp = \operatorname{closure}\left(\bigcup_{k \in \mathbb{N}} (\delta_k,\delta_{k+1})\right) \;.
\end{equation}
Each interval $(\delta_k,\delta_{k+1})$ is characterized by the property
that, for all $\delta \in (\delta_k,\delta_{k+1})$, $D_\delta^{L^\star}(s)$ has the same number of unstable
characteristic roots.
 
The boundaries $\delta_k$ of the stability intervals are called \emph{crossing points}. The method's key observation is that,
at crossing points, the quasipolynomial $D_{\delta_k}^{L^\star}(s)$ must have zeros at the imaginary axis. That is,
$D_{\delta_k}^{L^\star}(j\omega_\rc) = 0$ for some (not necessarily unique) $\omega_\rc \ge 0$, called a \emph{crossing 
frequency}. Note that, in the specific case of~\eqref{eq:D}, we have $D_{\delta_k}^{L^\star}(0) = l_n^\star \neq 0$, from 
which we conclude that $s = 0$ is not a root of $D_{\delta_k}^{L^\star}(s)$ and, consequently, $\omega_\rc > 0$. 
For strictly positive crossing frequencies, we can re-write $D_{\delta_k}^{L^\star}(j\omega_\rc) = 0$ as
\begin{equation} \label{eq:D_part}
    \eu^{j\omega_\rc \delta_k} = G(j\omega_\rc) \quad \text{with} \quad G(s) = -\frac{l_1^\star s^{n-1} + l_2^\star s^{n-2} + 
     \cdots + l_n^\star }{s^n} \;.
\end{equation}
Taking the modulus on both sides gives $1 = |G(j\omega_\rc)|$, which can be solved for $\omega_\rc$ independently of $\delta_k$. 
By taking the argument on both sides of~\eqref{eq:D_part}, we can easily find, for each crossing frequency, the crossing points
\begin{displaymath}
    \delta_k \in \left\{ \frac{2\pi i + \operatorname{Arg}(G(j\omega_\rc))}{\omega_\rc} \mid i \in \mathbb{N} \right\} \;, 
\end{displaymath}
where $\operatorname{Arg}(\cdot) \in [0,2\pi)$ returns the argument of a complex number.

\begin{figure}
\centering
\includegraphics[width=0.8\columnwidth]{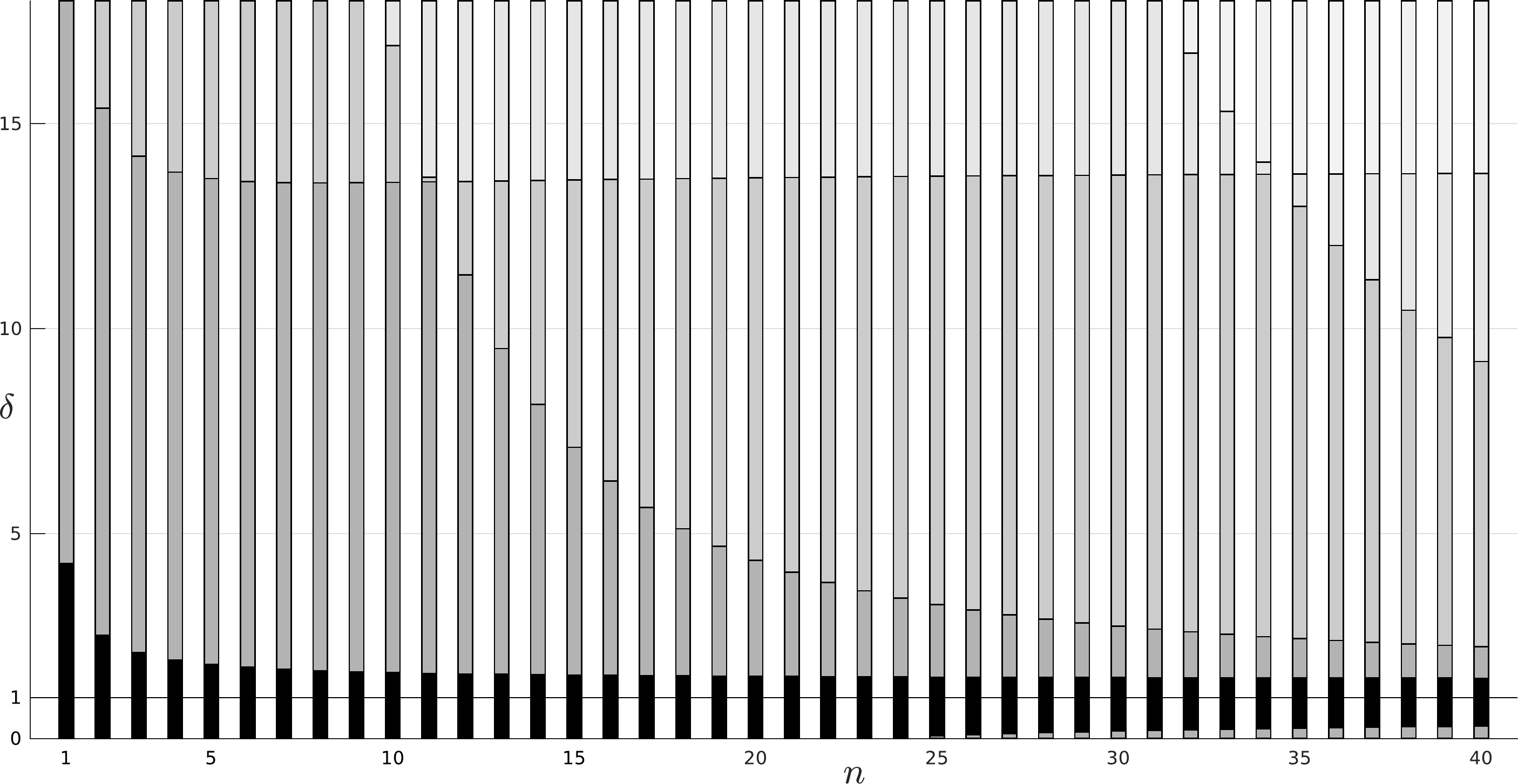}
\caption{Partition~\eqref{eq:partition} of the $\delta$-parameter space for the quasipolynomial $D_\delta^{L^\star}$~\eqref{eq:D}
 with $L^\star$ as in Theorem~\ref{thm:L}. The black intervals correspond to stable quasipolynomials.
 The number of unstable roots increases as the tone of the intervals decreases.}
\label{fig:d_vs_n}
\end{figure}

Recall that the \emph{crossing directions} are independent of the crossing frequencies \citep{cooke1986}. Thus, when $\omega_\rc$
is unique, the crossing direction is always towards instability (otherwise, for some finite $\delta$ there would be an infinite 
number of unstable roots, and this is impossible). In conclusion, the only open interval ensuring exponential stability
is the first one, $(0,\delta_1)$ with $\delta_1$ the \emph{delay margin} \citep{michiels}. Since $\delta = 1$ yields a stable
quasipolynomial, we know that $\delta_1 > 1$. This happens to be the case for the dimensions $1 \le n \le 8$.

For $9 \le n \le 25$ there are three crossing frequencies, one towards stability and two towards instability, whereas for
$26 \le n \le 46$, there are five crossing frequencies. At $n =23$, the polynomial $D_0^{L^\star}(s)$ ceases to be Hurwitz,
which implies that the first interval $(0,\delta_1)$ is no longer stable. It happens that the stable interval containing
$\delta = 1$ is the second one, as depicted in Figure~\ref{fig:d_vs_n}.

We see from Figure~\ref{fig:d_vs_n} that the predictor has reasonable robustness with respect to $\delta$ and hence with respect
to $\lambda$. However, the gain $\lambda$ cannot be set arbitrarily high as desired, and this compromises the
whole high-gain prediction scheme. Moreover, for high values of $n$ we need to make sure that $\delta$ (and hence $\lambda$) 
is not too low either.

\subsection{Cascading sub-predictors to recover the high-gain property}\label{sec:subpredictors}

Since the only limiting factor for the scalar gain is the delay, it is reasonable to substitute the single
predictor~\eqref{eq:high_gain_pred} by a chain of sequential sub-predictors \citep{germani2002,ahmed2012,najafi2013}. Indeed, when $N \in \mathbb{N}$
sequential sub-predictors are used, the \emph{effective delay} of each sub-predictor is reduced to
$h_\eu = h/N$ \citep{najafi2013}. A stable value in $\delta$-parameter space is now mapped to 
$\lambda = \delta/h_\eu = N \delta /h$. Thus, by choosing $N$ high enough, the scalar gain can be
increased without destroying the stability of $\psi_{1,\delta}^{L^\star}$. The idea is detailed in this subsection.

We define the sub-predictors
\begin{equation} \label{eq:sub_predictors}
\begin{aligned}
    \dot{\hat{x}}^1(t) &= A\hat{x}^1(t) + \Lambda_\lambda^r\left(L^\star\cdot\left(y(t)-C\hat{x}^1(t-h_\eu)\right)\right) +
     \varphi(\hat{x}^1(t),u(t-h+h_\eu)) \\
    \dot{\hat{x}}^j(t) &= A\hat{x}^j(t) + \Lambda_\lambda^r\left(L^\star C\cdot\left(\hat{x}^{j-1}(t) 
     - \hat{x}^j(t-h_\eu)\right)\right) + \varphi(\hat{x}^j(t),u(t - h + jh_\eu))
\end{aligned} \;,
\end{equation}
where $\hat{x}^1(t) \in \RE^n$, $\hat{x}^j(t) \in \RE^n$, and $j = 2,\dots,N$. To simplify the exposition we
define $\hat{x}^0(t) = x(t)$ and $e^0(t) = 0$ for all $t \ge -h$. It is not difficult to see that the errors
\begin{equation} \label{eq:sub_errors}
    e^j(t) = \hat{x}^{j}\left( t - jh_\eu \right) - \hat{x}^{j-1}\left( t - (j - 1)h_\eu \right) \;, \quad 
        j = 1,\dots,N \;, 
\end{equation}
evolve according to the time-delay differential equations
\begin{equation*} \label{eq:suberror_der}
    \dot{e}^j(t) = \psi_{\lambda,h_\eu}^{L^\star}(e^j_t) + \Delta\varphi\left(e^j_t(0),\hat{x}^{j-1}_t(-jh_\eu+h_\eu),u_t(-h)\right)
    - \Lambda_\lambda^r\left(LC e^{j-1}_t(-h_\eu)\right) \;.
\end{equation*}
For the dilated errors $\varepsilon^j(\tau) = \Lambda_{\lambda^{-1}}^r\left(e^j(\lambda^{-1}\tau)\right)$, $j=0,\dots,N$,
we have, similarly to Proposition~\ref{prop:high_gain_pred},
\begin{equation} \label{eq:finite_sub_error_delay}
    \frac{\rd\varepsilon^j(\tau)}{\rd \tau} = \psi_{1,\lambda h_\eu}^{L^\star}\left(\varepsilon^j_\tau\right) 
    + \kappa_\lambda^j\left(\varepsilon^j_\tau(0),\tau\right) + w_\lambda^j(\tau) \;, \quad j = 1,\ldots,N \;,
\end{equation}
where $\kappa_\lambda^j : \RE^n\times\REp \to \RE$ and $w_\lambda^j \in \LE(\REp,\RE^n)$ are defined by
\begin{displaymath}
    \kappa_\lambda^{j+1}\left(\xi,\tau\right) = 
     \lambda^{-1}\cdot\Delta\Phi_\lambda\left(\xi,\hat{x}^{j}_{\lambda^{-1}\tau}(-jh_\eu),u_{\lambda^{-1}\tau}(-h)\right) 
\end{displaymath}
$w_\lambda^{j+1}(\tau) = -LC \varepsilon^{j}_{\tau}(-\lambda h_\eu)$, $j=0,\ldots,N-1$.

We are now ready to state our main result.

\begin{thm}[Main result] \label{thm:main}
    Consider a uniformly observable system~\eqref{eq:normal_form} satisfying Assumption~\ref{ass:Lipschitz} and
    a chain of $N \in \mathbb{N}$ sub-predictors~\eqref{eq:sub_predictors} with $L^\star$ as in 
    Theorem~\ref{thm:L} and $r = [1,\ldots,n]$. Then, for $N$ large enough, there exists a gain $\lambda \ge 1$
	 such that
    \begin{equation} \label{eq:prediction_N}
        \lim_{t \to +\infty} \left( \hat{x}^{N}(t) - x(t+h) \right) = 0
    \end{equation}
    for any initial condition $x_0 \in \Omega$, initial errors $e^j_0 \in \CO([-h_\eu,0],\RE^n)$ 
	 with $e^j$ defined by~\eqref{eq:sub_errors}, and input $u \in \LE([-h,+\infty),U)$.
\end{thm}

\begin{proof}[Proof of the main result]
    Consider the dilated error equations~\eqref{eq:finite_sub_error_delay} and mark that 
    $\lim_{\tau \to +\infty} \varepsilon^j(\tau) = 0$ implies~\eqref{eq:prediction_N}.
    Since $\psi_{1,1}^{L^\star}$ is linear and exponentially stable, the virtual system
    \begin{displaymath}
         \frac{\rd\varepsilon(\tau)}{\rd \tau} = \psi_{1,1}^{L^\star}\left(\varepsilon_\tau\right) + w(\tau)
    \end{displaymath}
    is ISS with input $w \in \LE(\REp,\RE^n)$ and gain margin $\rho \in \CK_\infty$ (Corollaries~\ref{cor:linear}
    and~\ref{cor:robust}). On the other hand, it follows from Lemma~\ref{lem:Lip_homogeneous} that
    \begin{displaymath}
        \norm{\kappa^j_\lambda(\xi,\tau)} \le \frac{\gamma_\Phi}{\lambda} \norm{\xi}
    \end{displaymath}
    for all $\xi \in \RE^n$, $j = 1,\dots,N$, and $\lambda \ge 1$.
    By the linearity of $\psi_{1,1}^{L^\star}$, the gain margin can also be chosen to be linear.
    That is,
    \begin{equation} \label{eq:rho}
        \rho(\theta) = \gamma_m \cdot \theta
    \end{equation}
    for some $\gamma_m > 0$.
    Thus, finding the values of $N$ and $\lambda$ claimed in the theorem can be achieved in two simple steps:
    \begin{enumerate}
        \item Choose
         \begin{displaymath}
            \lambda^\star \ge \max\left\{ \frac{\gamma_\Phi}{\gamma_m}, 1 \right\}
         \end{displaymath}
         and note that the inequality implies that $\norm{\kappa^j_{\lambda}(\xi,\tau)} \le \rho(\norm{\xi})$ 
			for all $\lambda \ge \lambda^\star$.
        \item Choose $N \ge \lambda^\star h$ and set $\lambda = N/h$. 
    \end{enumerate}
    By Corollary~\ref{cor:robust}, the system
    \begin{displaymath}
     \frac{\rd\varepsilon(\tau)}{\rd \tau} = \psi_{1,1}^{L^\star}\left(\varepsilon_\tau\right) 
     + \kappa_\lambda^j\left(\varepsilon_\tau(0),\tau\right) + w(\tau)
    \end{displaymath}
    is ISS with input $w$, irrespective of the index $j$. Corollary~\ref{cor:cascade}
    ensures that the cascade~\eqref{eq:finite_sub_error_delay} is ISS with input $w^1$. Since $w^1(\tau) = 0$, we conclude
    that the trivial solution of~\eqref{eq:finite_sub_error_delay} is globally exponentially stable.
\end{proof}

\subsection{Numerical estimation of the gain margins} \label{sec:lmi}

Given an $n$-dimensional system~\eqref{eq:normal_form} with Lipschitz constant $\gamma_{\varphi}$, we can readily compute
the vector gain $L^\star$ as in Theorem~\ref{thm:L}. On the other hand, deciding on the number $N$ of required 
sub-predictors and the actual value of the scalar gain $\lambda$ (Steps 1 and 2 of the proof of Theorem~\ref{thm:main})
requires a gain margin~\eqref{eq:rho} for $\psi_{1,1}^{L^\star}$. 

An upper bound $\overline{\gamma}_m$ on $\gamma_m$ is $l^\star_n$, i.e.,  $\gamma_m \le \overline{\gamma}_m = l^\star_n$. 
Indeed, consider the perturbation 
\begin{displaymath}
 w(\tau) = 
    \begin{bmatrix}
        0 & \cdots & 0 & l^\star_n 
    \end{bmatrix}^\top C \varepsilon_\tau(0)
\end{displaymath}
and note that it satisfies the condition $\|w(\tau)\| \le l^\star_n \|\varepsilon_\tau(0)\|$. It is not difficult to see that
the delayed differential equation
\begin{displaymath}
    \frac{\rd\varepsilon(\tau)}{\rd \tau} = \psi_{1,1}^{L^\star}(\varepsilon_\tau) + w(\tau)
\end{displaymath}
has the characteristic equation $D_1^{L^\star}(s) - l_n^\star = 0$ , which has a solution at $s=0$.

\begin{figure}
\centering
\includegraphics[width=0.8\columnwidth]{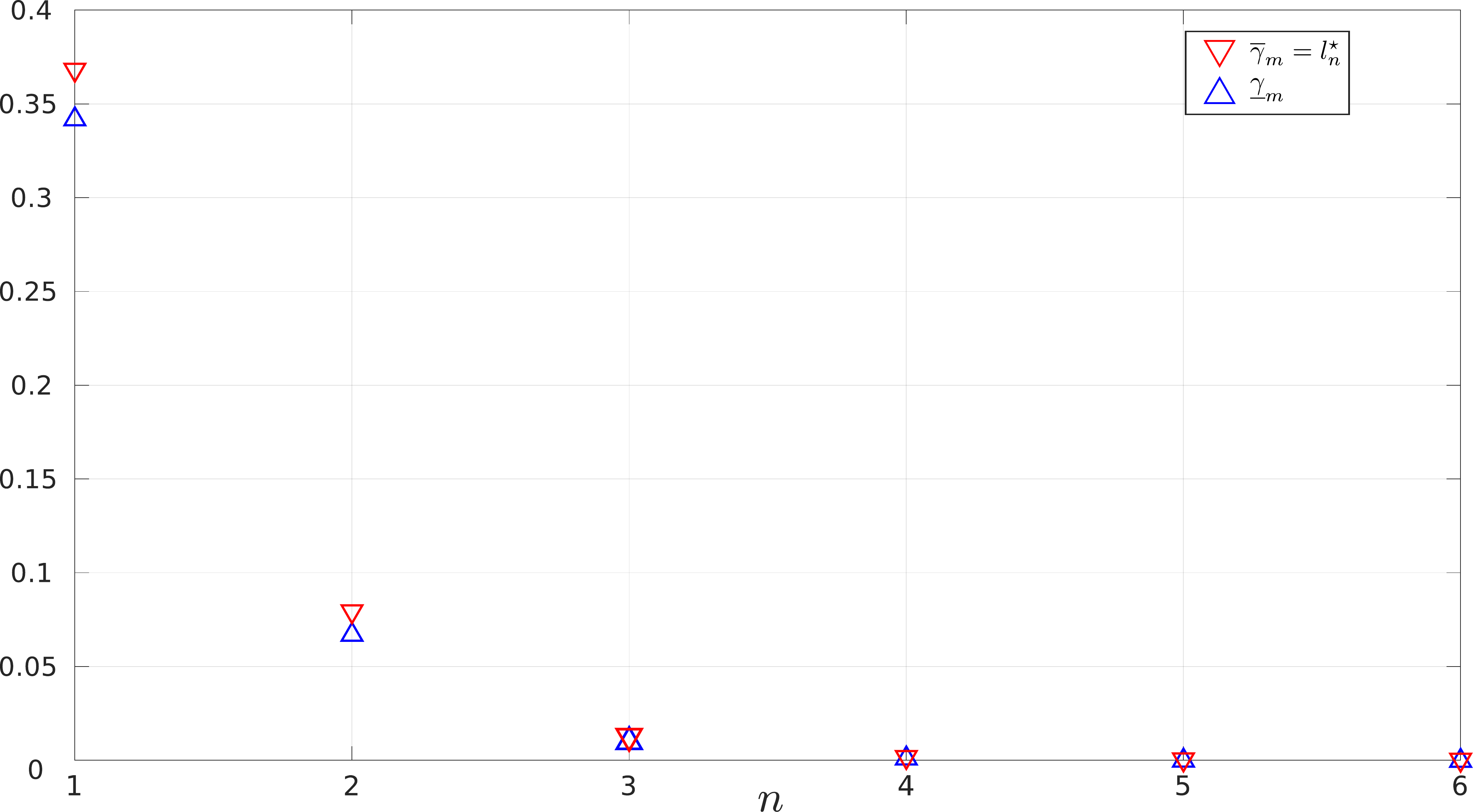}
\caption{Bounds on the gain margin of $\psi_{1,1}^{L^\star}$, defined by~\eqref{eq:psi}. The upper
 bound $\overline{\gamma}_m$ is given by $l^{\star}_n$ (downward triangle) and the lower bound $\underline{\gamma}_m$
 is obtained by Proposition~\ref{prop:margin} (upward triangle).}
\label{fig:gamma_vs_n}
\end{figure}

We capitalize on linear matrix inequality solvers to compute a lower bound on $\gamma_m$. 
We rely on Proposition~\ref{prop:margin}, which in turn relies on the LKF~\eqref{eq:LKF} given in the proof.
For each value of $n$, we maximize $\gamma_m$ subject to the existence of $P>0$, $R>0$, $S>0$, $P_2$, $P_3$, 
$P_4 \in \RE^{n\times n}$ such that $W < 0$ with $W$ as in~\eqref{eq:W}, $A_1 = -L^\star C$,
and $h=1$. Figure~\ref{fig:gamma_vs_n} shows the results obtained by means of the solver \cite{mosek}.

The upper and lower bounds are close, which certifies that the conservativeness introduced by the specific choice of~\eqref{eq:LKF}
is negligible. Unfortunately, we can see that $\gamma_m$ decreases with $n$. The latter attests to the difficulty that results from 
the compromise between the delay and gain margins, and to the fact that the compromise only worsens as the dimension of the problem
increases. We infer from Figure~\ref{fig:gamma_vs_n} that, for fixed $\gamma_\varphi$ and $h$, the number of required sub-predictors
increases exponentially with $n$. This is a possible explanation as to why all the examples we have found in the literature consider
predictors for systems of dimension no larger than three.

\section{Trade-offs between delay and gain margin} \label{sec:discussion}

A significant advantage of computing $L$ as described in Section~\ref{sec:L} is that 
the trade-off between $\gamma_\Phi$, $h$, $\lambda$, and $N$, can be formulated in very 
simple terms (cf. Steps 1 and 2 in the proof of Theorem~\ref{thm:main}). More importantly, there is no
conservatism introduced in the delay-margins obtained in Section~\ref{sec:delta}, and very little 
conservatism is introduced in the computation of the gain margins, as can be seen in Figure~\ref{fig:gamma_vs_n}.

As a reference for comparison, allow us to recall the trade-off derived by \cite{ahmed2012}. In our notation, 
they have
\begin{equation} \label{eq:ahmed}
\begin{aligned}
    \frac{h\lambda^3}{2} &> 2 \norm{P}\gamma_\Phi + h\lambda^2\big[\norm{A-LC}+2\norm{P}\gamma_\Phi\big]^2 \\
                       1 &> 2\norm{L}^2 \big[\norm{P}^2 + h^2\lambda^4\big]
\end{aligned} \;,  
\end{equation}
where $P$ is a symmetric positive-definite solution to the Lyapunov inequality
\begin{equation} \label{eq:lyap_ineq}
    P(A-LC) + (A-LC)^\top P \le -h\lambda^2 I_n
\end{equation}
(cf. inequalities~(36) in the last reference). 

\begin{prop} \label{prop:ahmed}
    Necessary conditions for the inequalities~\eqref{eq:ahmed},~\eqref{eq:lyap_ineq}, and $n \ge 2$ are
    \begin{equation} \label{eq:ahmed_2}
        h\lambda^2 < \frac{1}{\sqrt{2}n} \quad \text{and} \quad
        h \le \frac{1}{4\sqrt{2}n} \approx 0.1768\frac{1}{n} \;.
    \end{equation}
\end{prop}

Now we consider the trade-off outlined in~\cite{lei2016a}. For constant delays we have
\begin{equation} \label{eq:lei1}
     \sigma h \lambda \le 1 \;,
\end{equation}
where
\begin{equation} \label{eq:lei2}
    \sigma = \max\left\{ 8\norm{A - LC}^2\frac{\lambda_{\max}(P)}{\lambda_{\min}(P)}, 8\norm{PLC}^2, 2\sqrt{2}\norm{LC} \right\}
\end{equation}
and $P$ is the symmetric positive-definite solution to the Lyapunov equation 
\begin{displaymath}
    P(A-LC) + (A-LC)^\top P = - I_n
\end{displaymath}
(cf. inequality~(26) in the last reference).

\begin{prop} \label{prop:lei}
    A necessary condition for~\eqref{eq:lei1},~\eqref{eq:lei2}, and $n \ge 2$ is
    \begin{displaymath}
        h\lambda \le \frac{1}{8} \;.    
    \end{displaymath}
\end{prop}

These necessary (but not sufficient) conditions stand in contrast with the weaker, yet sufficient conditions
derived in this paper: $\lambda h = 1$ and $h$ free for $\gamma_\Phi = 0$. Also, it is not clear
from~\eqref{eq:ahmed}-\eqref{eq:lyap_ineq} nor from~\eqref{eq:lei1}-\eqref{eq:lei2} how to choose $L$ to optimize the
trade-off.

\section{Illustrative example} \label{sec:example}

\begin{figure}
    \centering
    \includegraphics[width=0.8\columnwidth]{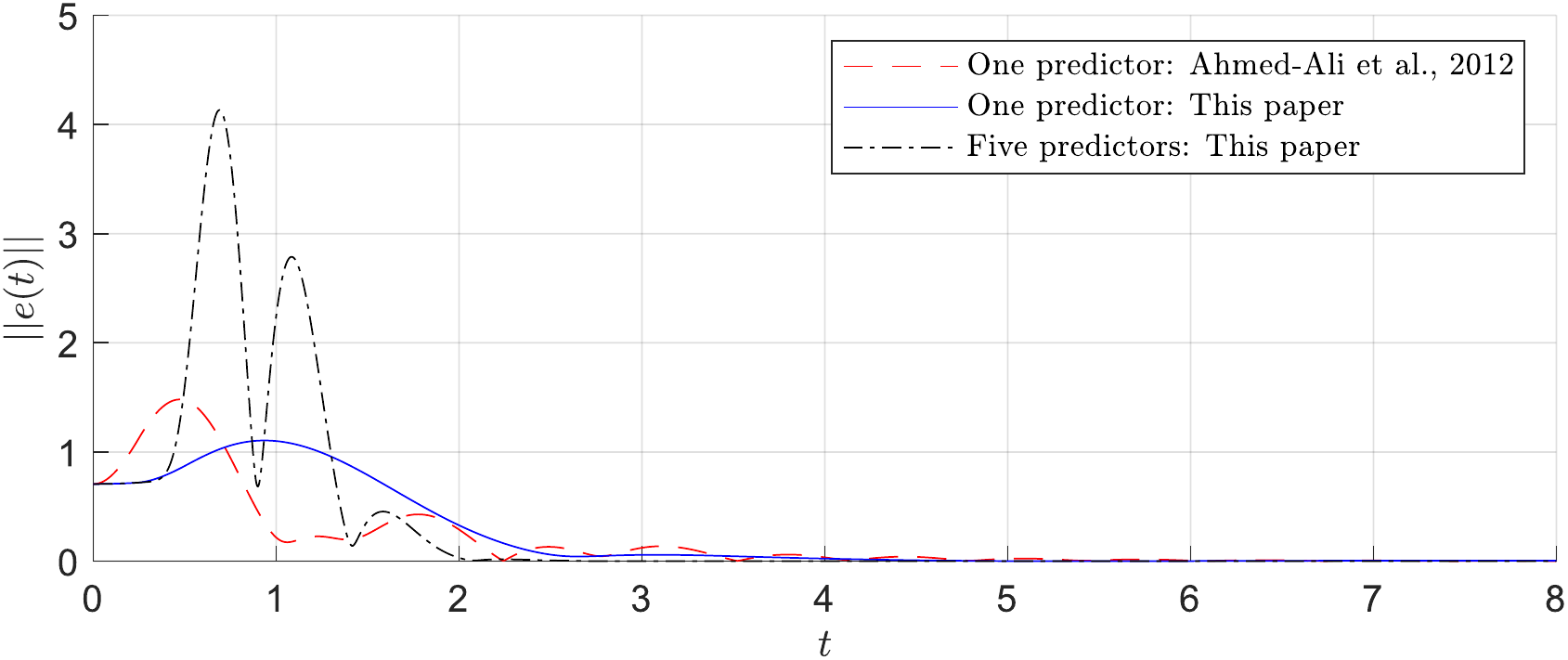}
    \caption{Estimation error $e(t)=\hat{x}^N(t-h)-x(t)$ for the chain of sub-predictors~\eqref{eq:sub_predictors}
     with $h=0.25$ using different tuning methods. Methodology proposed by \cite{ahmed2012} for $N=1$ predictor (dashed line).
     Methodology presented in this paper for $N=1$ (solid line) and $N=5$ (dot-dashed line).} 
    \label{fig:normError025}
\end{figure}

To illustrate the trade-offs discussed in Section~\ref{sec:discussion}, we revisit an example presented by \cite{ahmed2012}.
Transforming the original observation problem with delayed outputs into a prediction problem with delayed inputs gives the
system
\begin{equation} \label{eq:ahmed_example}
	\begin{aligned}
		\dot{x}_1(t) &= x_2(t) \\
		\dot{x}_2(t) &= -x_1(t) + 0.5\tanh(x_1(t)+x_2(t))+x_1(t)u(t-h) \;, 
	\end{aligned}  
\end{equation}
where the output is $y(t)=x_1(t)$. Note that~\eqref{eq:ahmed_example} is already in the canonical form~\eqref{eq:normal_form} with
\begin{equation*}
    \varphi(x,u) = \begin{bmatrix}
        0 & -x_1 + 0.5\tanh(x_1+x_2)+x_1 u
    \end{bmatrix}^\top \;,
\end{equation*}
which clearly satisfies Assumption~\ref{ass:Lipschitz}. A Lipschitz constant $\gamma_\Phi$
can be computed as described in the proof of Lemma~\ref{lem:Lip_homogeneous}.
Note that $\gamma_1 = 0$ and $\gamma_2$ is an upper bound for
\begin{equation*}
    \norm{\frac{\partial \varphi_2}{\partial x}} = \sqrt{1 - 2 u + u^2 +(u-1)
     \operatorname{sech}^2(x_1+x_2) + 0.5 \operatorname{sech}^4(x_1+x_2)} \;.
\end{equation*}
Considering the input proposed by \cite{ahmed2012}, $u(t) = 0.1 \sin(0.1t)$, which satisfies the bound $|u(t)| \leq 0.1$ for all $t \in \RE$, and noting that $0 < \operatorname{sech}(x_1+x_2) \leq 1$ for all $x_1, x_2 \in \RE$, a Lipschitz constant~\eqref{eq:gamma_Phi} is
\begin{equation*}
    \gamma_\Phi = \gamma_2 = \sqrt{1 + 2|u| + u^2} = \sqrt{1.21} \;.
\end{equation*}
\cite{ahmed2009,ahmed2012} propose the use of one predictor with the gains $\lambda = 2$ and $L=\begin{bmatrix}2 & 1 \end{bmatrix}^\top$. Having two eigenvalues at $-1$, the matrix $(A-LC)$ is Hurwitz. However, we remark that conditions~\eqref{eq:ahmed} are not satisfied (indeed, the conditions cannot hold for any $L$, as~\eqref{eq:ahmed_2} is transgressed for $h=0.25$). Thus, more than one sub-predictor is required to satisfy their proposed sufficient conditions. Nevertheless, simulations show the convergence of the estimation errors to zero (Figure~\ref{fig:normError025}).

\begin{figure}
    \centering
    \includegraphics[width=0.6\columnwidth]{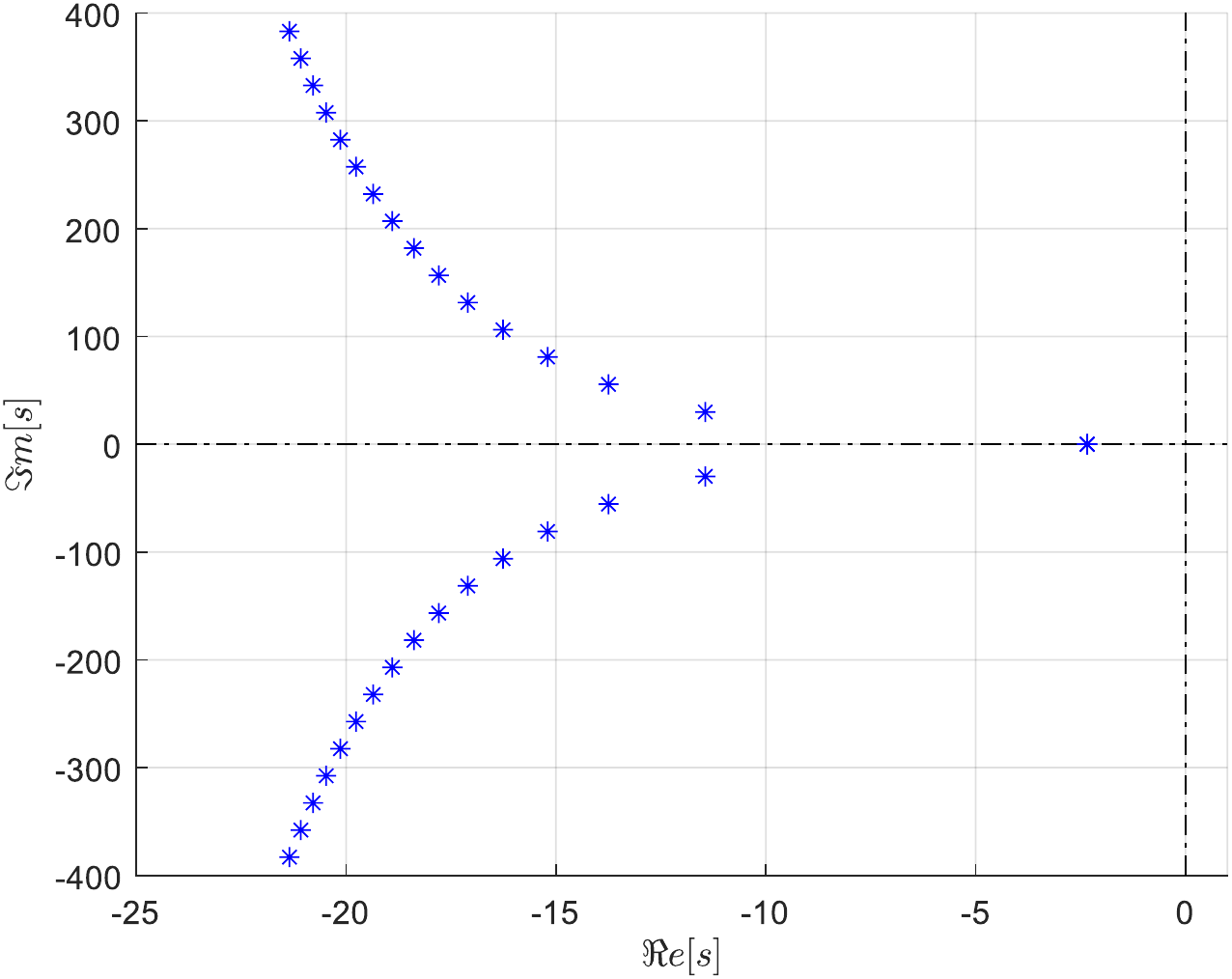}
    \caption{Spectrum of $\psi_{\lambda,h}^{L^\star}\left(e_t\right)$ for $n=2$, $\lambda=1/h$ and $h=0.25$. 
     The functional vector field is stable, as there is a dominant root of multiplicity $3$ at
     $\sigma_\star/h < 0$.}
    \label{fig:spectrum025}
\end{figure}

We now tune the sub-predictors~\eqref{eq:sub_predictors}
with $L^\star=\begin{bmatrix}l_1^\star & l_2^\star\end{bmatrix}^\top$ as in Theorem~\ref{thm:L}.
The rightmost root of \eqref{eq:q} is  
$\sigma_\star=-2+\sqrt{2}$, which gives $l_1^\star=0.4612$ and $l_2^\star=0.0791$.
Recall that this gain assigns a dominant root with multiplicity $n+1=3$ at $\sigma_\star/h$ to the spectrum of $\psi_{1/h,h}^{L^\star}\left(e_t\right)$.
The spectrum was computed using QPmR \citep{vyhlidal2014} and is shown in Figure~\ref{fig:spectrum025}.

As discussed in Section~\ref{sec:subpredictors}, there is a simple way to establish a number of sub-predictors guaranteeing the exponential convergence of the prediction error. For $n=2$, we have $\gamma_m=0.0673$ (Figure~\ref{fig:gamma_vs_n}), so $\lambda^\star = \max\left\{ \gamma_\Phi/\gamma_m,\,1\right\} = 16.3556$. 
For the delay $h=0.25$ we have $\lambda^\star h = 4.0889$, from where we conclude that $N = 5$ sub-predictors with the scalar gain $\lambda = N/h = 20$ are sufficient to solve Problem~\ref{prob:prediction}. It is worth mentioning that, for $N = 5$, there is no $\lambda$ satisfying~\eqref{eq:ahmed_inter}, so more than five sub-predictors are required to satisfy the conditions derived by \cite{ahmed2012}.

\begin{figure}
    \centering
    \includegraphics[width=0.8\columnwidth]{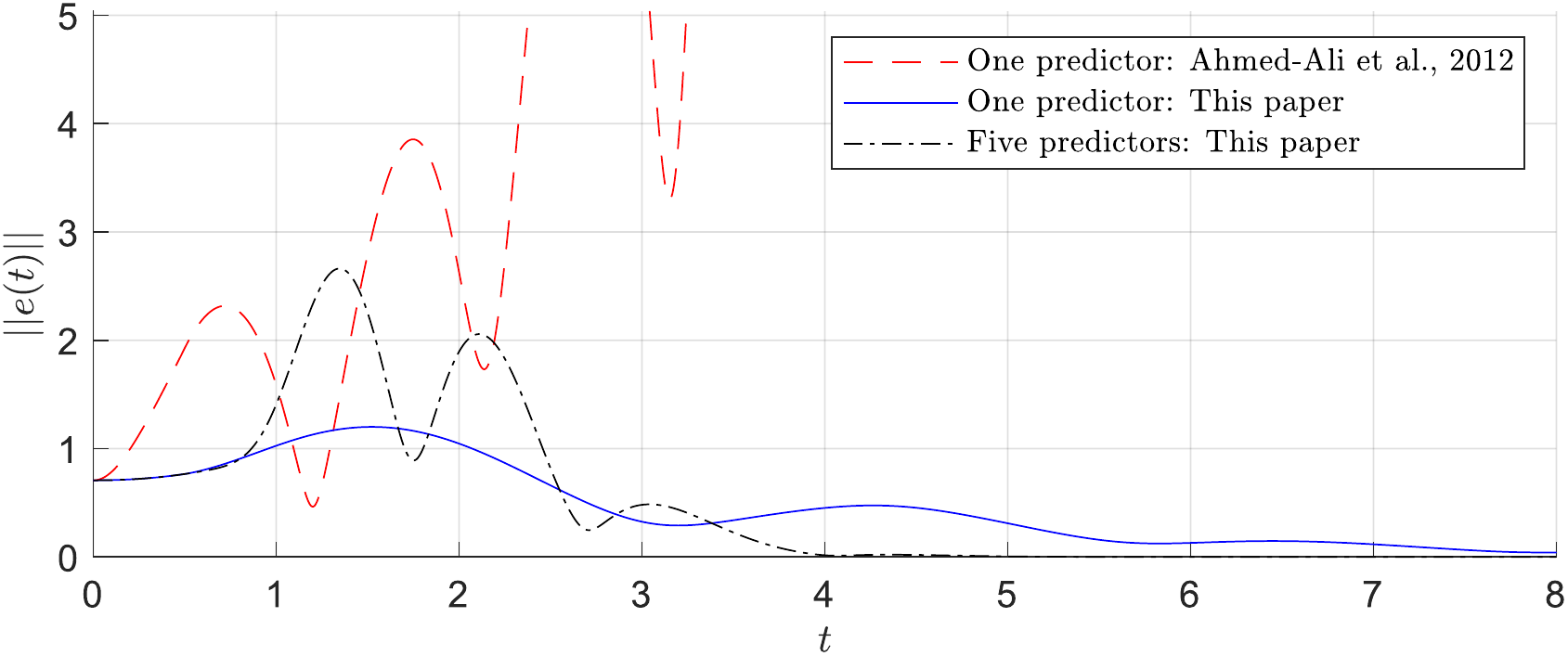}
    \caption{Estimation error $e(t)=\hat{x}^N(t-h)-x(t)$ for the chain of sub-predictors~\eqref{eq:sub_predictors}
     with $h=0.5$ using different tuning methods. Methodology proposed by \cite{ahmed2012} for $N=1$ predictor (dashed line).
     Methodology presented in this paper for $N=1$ (solid line) and $N=5$ (dot-dashed line).} 
    \label{fig:normError050}
\end{figure}

For comparison purposes, we simulate a chain of sub-predictors with the tuning rules just described, both for $N =5$, $\lambda = N/h = 20$ and for $N = 1$, $\lambda = N/h = 4$.  These are presented in Figure~\ref{fig:normError025}. When using five cascaded sub-predictors, the error convergence rate is faster, but it results in larger oscillations. Although a predictor alone does not satisfy our sufficient conditions either, the simulations also show
convergence to zero with a fast convergence rate but with much better damping. To test the system's resilience with respect to changes in the delay, we repeat the three simulations described above, but with a delay $h = 0.5$. The results are presented in Figure~\ref{fig:normError050}. The error dynamics for the predictor tunned using the methodology presented by \cite{ahmed2012} are unstable, whereas the tunning presented in
this paper yields stable dynamics, both for $N = 1$ and $N=5$.

 

\section{Conclusions}


We introduce a new method for tuning high-gain observer--predictors for uniformly observable nonlinear systems with input delay. In contrast with previous proposals, which design the observer for the delay-free system and treat the delay as a disturbance, we use the multiplicity-induced-dominance pole placement to tune the gain. Nonlinearities with large Lipschitz constants and large delays may impose a cascade of sub-predictors. The input-to-state stability of the nonlinear time-delay error dynamics is analyzed by Lyapunov-Krasovskii techniques.

The proposed approach is shown to provide a less conservative trade-off between delay and gain margin. The explicit calculation of the gain margin leads to a simple, yet non-conservative criterion for determining the required number of sub-predictors. 

It is worth mentioning that, for time-delay systems, a downside of the nonlinear domination approach is that the gain margin decreases exponentially with the system's dimension, which may result in an unreasonably large number of sub-predictors.

\section*{Acknowledgment} 
To Prof. 
C. F. M\'endez-Barrios for his invaluable advice on the stability of quasipolynomials.

\appendix

\section{Proofs}

\begin{proof}[Proof of Lem.~\ref{lem:homo_linear}]
    For $A$ a Jordan block with zero eigenvalue we have
    \begin{displaymath}
        \left( A\Lambda_\lambda^r(x) \right)_i = \lambda^{r_{i+1}}\cdot x_{i+1} = \lambda\cdot\lambda^{r_i}\cdot x_{i+1} \quad \text{and} \quad
        \left(\Lambda_\lambda^r(A x ) \right)_i = \lambda^{r_i} \cdot x_{i+1} \;, \quad i = 1,\dots,n-1 \;. 
    \end{displaymath}
    Also,
    \begin{displaymath}
        \left( A\Lambda_\lambda^r(x) \right)_n = \left( \Lambda_\lambda^r(A\cdot x) \right)_n = 0 \;,
    \end{displaymath}
    so $A\Lambda_\lambda^r(x) = \lambda\cdot\Lambda_\lambda^r(Ax)$. By direct computation, $C\Lambda_\lambda^r(x) = \lambda \cdot C x = \lambda \cdot x_1$.
\end{proof}

\begin{proof}[Proof of Lem.~\ref{lem:Lip_homogeneous}]
    For $x \in \RE^n$ and $1 \le i \le n$ we introduce the notation $x_{[1,i]} = [x_1,\dots,x_i]$. Note that
    $\varphi$ is globally Lipschitz in $x \in \Omega$, uniformly in $u \in U$ if, and only if, each of its
    components is uniformly globally Lipschitz, that is, if 
    $|\varphi_i(\hat{x}_{[1,i]},u) - \varphi_i(_{[1,i]},u)| \le \gamma_i \norm{\hat{x}_{[1,i]} - x_{[1,i]}}$
    for some nonnegative $\gamma_i$, $i = 1,\dots,n$. Thus,
    \begin{displaymath}
        \left|\Delta\varphi_i\left(\Lambda_\lambda^{r_{[1,i]}}(e_{[1,i]}),x_{[1,i]},u\right)\right| 
         \le \gamma_i \norm{\Lambda_\lambda^{r_{[1,i]}}(e_{[1,i]})} \;.
    \end{displaymath}
    On the other hand, for these particular weights we have
    \begin{displaymath}
        \norm{\Lambda_\lambda^{r_{[1,i]}}(e_{[1,i]})} \le \lambda^{r_i}\norm{e_{[1,i]}} \;, \quad \lambda \ge 1 \;,
    \end{displaymath}
    so that 
    \begin{equation*}
        \left|\Delta\varphi_i\left(\Lambda_\lambda^{r_{[1,i]}}(e_{[1,i]}),x_{[1,i]},u\right)\right| 
         \le \gamma_i \cdot \lambda^{r_i}\norm{e_{[1,i]}} \;.
    \end{equation*}
    This implies~\eqref{eq:lip_l} with 
    \begin{equation}\label{eq:gamma_Phi}
        \gamma_{\Phi} = \sqrt{\sum_{i=1}^n \gamma_i^2}.
    \end{equation}
\end{proof}

\begin{proof}[Proof of Lem.~\ref{lem:pseudo_homo}]
    Direct use of~\eqref{eq:psi} gives
    \begin{displaymath}
        \psi_{\lambda,\bar{h}}^L\left( \Lambda_\lambda^r\circ e_t \right) = A\Lambda_\lambda^r\left( e_t(0) \right) - 
         \Lambda_\lambda^r\left( LC \Lambda_\lambda^r(e_t(-\bar{h})\right) \;.
    \end{displaymath}
    By Lemma~\ref{lem:homo_linear} we have
    \begin{displaymath}
        \psi_{\lambda,\bar{h}}^L\left( \Lambda_\lambda^r\circ e_t \right) = \lambda\cdot \Lambda_\lambda^r\left( A e_t(0) \right) - 
         \Lambda_\lambda^r\left( LC (\lambda\cdot e_t(-\bar{h}))\right)
    \end{displaymath}
    and, by the linearity of $\Lambda_\lambda^r$, we recover~\eqref{eq:pseudo_homo}.
\end{proof}

\begin{proof}[Proof of Prop.~\ref{prop:dilation}] 
    We have $\dot{\xi}(t) = \lambda^d \cdot \Lambda_\lambda^r \left( f\left( x(\lambda^{d}\cdot t)\right) \right) = 
    f\left( \Lambda_\lambda^r\left( x(\lambda^{d}\cdot t)\right) \right) = f\left(\xi(t)\right)$.
\end{proof}

\begin{proof}[Proof of Prop.~\ref{prop:margin}]
    It suffices to show that the origin of~\eqref{eq:general} with
    \begin{equation} \label{eq:linear_pert}
        f(e_t,w) = Ae_t(0) + A_1e_t(-h) + w
    \end{equation}
    and $w(t) = \kappa(e_t(0),t)$
    is asymptotically stable for any $\kappa:\RE^n\times\RE \to \RE^n$ such that $\norm{\kappa(\xi,t)} \le \gamma_m\norm{\xi}$. 
    Consider the LKF \citep{fridmanE}
    \begin{equation} \label{eq:LKF}
        V\left( \phi \right) = \phi(0)^\top P\phi(0)+
        \int_{-h}^{0} \phi(\theta)^\top S\phi(\theta)\rd \theta + 
        h\int_{-h}^{0} \int_{\vartheta}^{0} \frac{\rd \phi(\theta)}{\rd \theta}^\top R \frac{\rd \phi(\theta)}{\rd \theta} 
        \rd \theta\rd \vartheta \;.
    \end{equation}
    Driver's derivative in the direction $v \in \RE^n$ is
    \begin{multline*}
        D^+V(\phi,v) = 2\phi(0)^\top P v + \phi(0)^\top S\phi(0) + h^2 v^\top R v - \phi(-h)^\top S\phi(-h) \\
        - h \int_{-h}^{0}\frac{\rd \phi(\vartheta)}{\rd \theta}^\top R \frac{\rd \phi(\vartheta)}{\rd \theta}
        \rd \vartheta \;.
    \end{multline*}
    Application of Jensen's inequality to the integral term gives the bound
    \begin{multline*}
        D^+V(\phi,v) \le 2\phi(0)^\top P v + \phi(0)^\top S\phi(0) + h^2 v^\top R v - \phi(-h)^\top S\phi(-h) \\
        - \left( \phi(0) - \phi(-h) \right)^\top R \left( \phi(0) - \phi(-h) \right) \;.
    \end{multline*}
    Thus, the Dini derivative of $V(e_t)$ along the trajectories of the vector field $f(e_t,w)$ satisfies
    \begin{multline} \label{eq:dini}
        D^+V(e_t,f(e_t,w(t))) \le 
         2e_t(0)^\top P v + e_t(0)^\top Se_t(0) + h^2 v^\top(t) R v(t) - e_t(-h)^\top Se_t(-h) \\
        - \left( e_t(0) - e_t(-h) \right)^\top R \left( e_t(0) - e_t(-h) \right) \Big|_{v(t) = f(e_t,w(t))}
    \end{multline}
    almost everywhere. Instead of making the substitution $v(t) = f(e_t,w(t))$, we follow the descriptor method
    \citep{fridmanE} and propose the null term
    \begin{equation} \label{eq:descriptor}
        2\left[ P_2e_t(0) + P_3v(t) + P_4 w(t)\right]^\top \left[ f(e_t,w(t)) - v(t) \right] = 0
    \end{equation}
    with $f$ given by~\eqref{eq:linear_pert}. By adding~\eqref{eq:descriptor} and the gain margin hypothesis
    $\gamma_m^2\norm{e_t(0)}^2 - \norm{w(t)}^2 \ge 0$ to~\eqref{eq:dini} we obtain the bound
    \begin{displaymath}
        D^+V(e_t,f(e_t,w(t))) \le E(t)^\top W E(t)
    \end{displaymath}
    with $W$ as in~\eqref{eq:W} and $E(t)^\top = \begin{bmatrix} e_t(0)^\top & v(t)^\top & e_t(-h)^\top & w(t)^\top \end{bmatrix}$.
    We conclude that the hypothesis $\norm{\kappa(e_t(0),t)} \le \gamma_m\norm{e_t(0)}$ implies
    the asymptotic stability of $\dot{e}(t) = Ae_t(0)+A_1e_t(-h)+\kappa(e_t(0),t)$ whenever $W < 0$. 
\end{proof}

\begin{proof}[Proof of Prop.~\ref{prop:high_gain_pred}]
    Differentiating~\eqref{eq:dilated_error} with respect to $\tau$ gives
    \begin{equation} \label{eq:dilated_error_der}
        \frac{\rd\varepsilon(\tau)}{\rd \tau} = \lambda^{-1}\cdot\Lambda_{\lambda^{-1}}^r\left(\dot{e}(\lambda^{-1}\tau)\right) \;,
    \end{equation}
    We can readily see that $\psi_{\lambda,\bar{h}}^L\left(e_{\lambda^{-1}\tau}\right) = 
    \psi_{\lambda,\lambda \bar{h}}^L\left( \Lambda_\lambda^r\circ\varepsilon_\tau \right)$ so, by
    Lemma~\ref{lem:pseudo_homo} we have
    \begin{equation} \label{eq:psi_epsilon}
        \psi_{\lambda,\bar{h}}^L\left(e_{\lambda^{-1}\tau}\right) = 
         \lambda\cdot\Lambda_{\lambda}^r\left(\psi_{1,\lambda \bar{h}}^L(\epsilon_\tau)\right) \;.
    \end{equation}
    Substituting~\eqref{eq:psi_epsilon} in~\eqref{eq:error_der} gives
    \begin{equation*}
        \dot{e}(\lambda^{-1} \tau) =  \lambda\cdot\Lambda_{\lambda}^r\left(\psi_{1,\lambda \bar{h}}^L(\epsilon_\tau)\right) + 
         \Delta\varphi\left(\Lambda_\lambda^r(\varepsilon(\tau)),x(\lambda^{-1}\tau),u(\lambda^{-1}\tau-h)\right) \;,
    \end{equation*}
    so that
    \begin{equation*}
        \frac{\rd\varepsilon(\tau)}{\rd \tau} = \psi_{1,\lambda \bar{h}}^L(\epsilon_\tau) + 
         \lambda^{-1}\Lambda_{\lambda^{-1}}^r\left(
          \Delta\varphi\left(\Lambda_\lambda^r(\varepsilon(\tau)),x(\lambda^{-1}\tau),u(\lambda^{-1}\tau-h)\right)
         \right) \;.
    \end{equation*}
\end{proof}

\begin{proof}[Proof of Prop.~\ref{prop:ahmed}]
    The Lyapunov inequality~\eqref{eq:lyap_ineq} establishes the bound $2\norm{P} \ge nh\lambda^2/l_1$ \citep{patel1978}.
    Also, since
    \begin{displaymath}
        (A-LC)^\top(A-LC) = 
        \begin{bmatrix}
            \norm{L}^2 &    -l_1 & \cdots & -l_{n-1} \\
                  -l_1 &       1 & \cdots & 0 \\
                \vdots &  \vdots & \ddots & \vdots \\
                -l_{n-1} &       0 & \cdots & 1
        \end{bmatrix} \;,
    \end{displaymath}
    we have 
    \begin{equation} \label{eq:A_LC}
        \norm{A-LC} \ge \max\left\{ \norm{L}, 1 \right\} \;.
    \end{equation}
    Thus, the inequalities~\eqref{eq:ahmed} imply
    \begin{subequations} \label{eq:ahmed_inter}
    \begin{align}
        \frac{\lambda}{2} &> \frac{n}{l_1}\gamma_\Phi + \left[\max\left\{ \norm{L}, 1 \right\} 
         +\frac{nh\lambda^2}{l_1}\gamma_\Phi\right]^2 
          \label{eq:ahmed_inter_1} \\
                        1 &> 2\frac{\norm{L}^2}{l_1^2} h^2\lambda^4 \left[n^2 + l_1^2\right] \;. \label{eq:ahmed_inter_2}
    \end{align}
    \end{subequations}
    Since $\norm{L}^2/l_1^2 \ge 1$,~\eqref{eq:ahmed_inter_2} implies the 
    first inequality in~\eqref{eq:ahmed_2}. By setting $\gamma_{\Phi} = 0$ in~\eqref{eq:ahmed_inter_1} we obtain
    $\lambda > 2$, which directly implies the second bound in~\eqref{eq:ahmed_2}.
\end{proof}

\begin{proof}[Proof of Prop.~\ref{prop:lei}]
    Using~\eqref{eq:A_LC} we find the simple bound $\sigma \ge 8$ for~\eqref{eq:lei2}. The claim follows directly
    from~\eqref{eq:lei1}.
\end{proof}

\begin{proof}[Proof of Cor.~\ref{cor:high_gain_obs}]
    Set $\bar{h} = 0$ in Proposition~\ref{prop:high_gain_pred} and note that $\Upsilon_\lambda^L(e(t)) = \psi_{\lambda,0}^L(e_t)$.
\end{proof}

\begin{proof}[Proof of Cor.~\ref{cor:robust}]
    By Theorem~\ref{thm:ISS}, system~\eqref{eq:general} admits a LKF satisfying~\eqref{eq:LKF_ISS}.
    Suppose that $\chi \in \CK_\infty$ (otherwise, replace $\chi(\theta)$ by $\chi(\theta)+\theta$).
    For $\theta \in \REp$ set
    \begin{displaymath}
        \rho(\theta) = \frac{1}{2}\cdot\chi^{-1}\left(\frac{1}{2}\cdot \alpha_1(\theta) \right) \;, \quad
        \tilde{\chi}(\theta) = 2 \cdot\chi\left( 2 \theta \right) \;,
    \end{displaymath}
    and note that $\rho, \tilde{\chi} \in \CK_\infty$. By the triangle and the weak triangle inequalities \citep{kellet2014},
    we have
    \begin{equation} \label{eq:triangle}
        \chi(\norm{\kappa(\phi(0),t) + v}) \le \chi(\norm{\kappa(\phi(0),t)} + \norm{v}) \le 
            \chi(2\norm{\kappa(\phi(0),t)}) + \chi(2\norm{v})
    \end{equation}
    for all $\phi(0), v \in \RE^n$. On the other hand, the inequalities
    \begin{equation} \label{eq:hypot}
        \norm{\kappa(\phi(0),t)} \le \rho(\norm{\phi(0)}) \quad \text{and} \quad 
         \tilde{\chi}(\norm{v}) \le V(\phi)
    \end{equation}
    are equivalent, respectively, to
    \begin{displaymath}
        \chi\left( 2 \norm{\kappa(\phi(0),t)}\right) \le \frac{1}{2}\alpha_1(\norm{\phi(0)}) \quad \text{and} \quad 
         \chi(2\norm{v}) \le \frac{1}{2}V(\phi) \;.
    \end{displaymath}
    Thus, using~\eqref{eq:triangle} we can see that~\eqref{eq:hypot} implies $\chi(\norm{\kappa(\phi(0),t) + v}) \le V(\phi)$.
    From~\eqref{eq:LKF_ISS} we have
    \begin{displaymath}
        D^+ V(\phi,f(\phi,\kappa(\phi(0),t)+v)) \le -\alpha_3(\norm{\phi(0)}) \;.
    \end{displaymath}
    Invoking Theorem~\ref{thm:ISS} again, we see that~\eqref{eq:feed_pert} is also ISS (although with a higher ISS gain
    since $\tilde{\chi}(\theta) > \chi(\theta)$ for $\theta > 0$).
\end{proof}


\bibliography{biblio}

\end{document}